\documentclass[12pt,english]{article}
\usepackage{theorem}
\usepackage{babel}
\usepackage{mathrsfs}
\usepackage{a4wide}
\usepackage{graphics}
\usepackage{epsfig}
\usepackage{amssymb}
\usepackage{latexsym}
\usepackage{amsmath}

\newtheorem{theorem}{Theorem}
\newtheorem{proposition}[theorem]{Proposition}
\newtheorem{lemma}[theorem]{Lemma}

\title{Surjectivity of quotient maps \\ for
algebraic $\CP$-actions \\ and polynomial maps \\ with contractible fibres}
\author{Philippe Bonnet}
\date{}

\usepackage{amsmath}
\let\cal\mathcal
\newcommand{\dem}{{\em Proof: }}
\newcommand{\qed}{\begin{flushright} $\blacksquare$\end{flushright}}

\newcommand{\der}{ \partial}
\newcommand{\CC}{\mathbb C}

\newcommand{\CP}{(\mathbb{C} , +)}

\newcommand{\CX}{\mathbb{C} [x_1,..,x_n]}

\begin{document}

\maketitle

\begin{center} { \small
Departamento de Algebra,\\
 Geometria y Topologia \\
Universidad de Valladolid \\
47005 Valladolid, Spain. \\
e-mail: pbonnet@agt.uva.es}
\end{center}

\begin{abstract}
In this paper, we establish two results concerning algebraic
$\CP$-actions on $\CC^n$. First of all, let $\varphi$ be an
algebraic $(\CC,+)$-action on $\CC^3$. By a result of Miyanishi,
its ring of invariants is isomorphic to $\CC[t_1,t_2]$. If
$f_1,f_2$ generate this ring, the quotient map of $\varphi$ is the
map $F:\CC^3 \rightarrow \CC ^2, \; x\mapsto (f_1(x),f_2(x))$. By
using some topological arguments, we prove that $F$ is always
surjective. Secondly we are interested in dominant polynomial maps
$F:\CC^n \rightarrow \CC^{n-1}$ whose connected components of
their generic fibres are contractible. For such maps, we prove the
existence of an algebraic $\CP$-action $\varphi$ on $\CC^n$ for
which $F$ is invariant. Moreover we give some conditions so that
$F^*(\CC[t_1,..,t_{n-1}])$ is the ring of invariants of $\varphi$.
\end{abstract}

\section{Introduction}

In this paper, we are going to study some properties of algebraic
$\CP$-actions on $\CC^n$. An algebraic $(\CC,+)$-action on $\CC^n$
is a regular map $\varphi: \CC \times \CC^n \rightarrow \CC^n$
such that $\varphi(u;\varphi(v;x))=\varphi(u+v;x)$ for all
$u,v,x$. It is well-known that $\varphi$ is obtained by
integrating a locally nilpotent derivation $\der$ on $\CX$, that
is a derivation $\der$ such that, for any polynomial $R$, there
exists an integer $k>0$ such that $\der ^k (R)=0$. A polynomial
$R$ is invariant if $R\circ \varphi=R$, or equivalently if
$\der(R)=0$. These polynomials form a ring called the ring of
invariants of $\varphi$, and denoted by $\CX ^{\varphi}$. We say
that $\varphi$ {\em satisfies condition $(H)$} if its ring of
invariants is isomorphic to a polynomial ring in $(n-1)$
variables. In this case, $\varphi$ is provided with a quotient map
$F$ defined as follows: If $f_1,..,f_{n-1}$ denote a system of
generators of $\CX ^{\varphi}$, then $F$ is the map: $$ F:\CC^n
\longrightarrow \CC ^{n-1}, \; x\longmapsto (f_1(x),..,f_{n-1}(x))
$$ Note that for $n>3$, the assumption $(H)$ need not be satisfied
([Wi]). Conversely a dominant polynomial map $F=(f_1,..,f_{n-1})$
is the {\em quotient map of a $\CP$-action on $\CC^n$} if there
exists an algebraic $\CP$-action $\varphi$ on $\CC^n$ such that:
$$\CC[f_1,..,f_{n-1}]= \CX ^{\varphi}$$ First of all, we establish
a property concerning algebraic $\CP$-actions on $\CC^3$.
According to a result of Miyanishi (see~\cite{Miy}), such an
action always satisfies condition $(H)$, and is therefore provided
with a quotient map $ F:\CC^3 \rightarrow \CC ^2$. In \cite{Kr},
Kraft conjectures that every fixed-point free $(\CC,+)$-action
$\varphi$ on $\CC^3$ is trivial, which means that it is conjugate
via an automorphism of $\CC^3$ to the action: $$
\varphi^0(t;x_1,x_2,x_3)=(x_1 + t,x_2,x_3) $$ This is known if its
quotient space is separated. More precisely, a fixed-point free
$(\CC,+)$-action $\varphi$ on $\CC^3$ is trivial if and only if:

\begin{itemize}
\item{$F$ is non-singular,}
\item{$F$ is surjective,}
\item{Every fibre of $F$ is connected.}
\end{itemize}
Daigle proved in \cite{Da} that $F$ is {\em non-singular in
codimension 1}, i.e. its singular set has codimension $\geq 2$.
Moreover he derived a jacobian formula for the locally nilpotent
derivation generating $\varphi$. His formula implies in particular
that $F$ is non-singular if $\varphi$ is fixed-point free. In an
attempt to understand the behaviour of $(\CC,+)$-actions on
$\CC^3$, we are going to study the second condition on $F$ given
above. More precisely:

\begin{theorem} \label{surj}
Let $\varphi$ be any algebraic $(\CC,+)$-action on $\CC^3$. Then its
quotient map $F$ is surjective.
\end{theorem}
Consequently a fixed-point free $(\CC,+)$-action $\varphi$ on $\CC^3$ is
trivial
if and only if every fibre of $F$ is connected.

The proof uses both algebraic and topological methods. First we
check that the complement of the image of $F$ is at most finite.
We assume that this complement is not empty, and denote by $x$ one
of its points. Let $K$ be an homological 3-sphere, i.e. a singular
3-cycle whose class generates the group $H_3(\CC^2 -\{x\})$. We
construct a singular 3-cycle $\Sigma$ in $\CC^3$, such that $F$
maps its homological class $[\Sigma]$ in $\CC^3$ $p$ times on
$[K]$, where $p$ is a positive integer. In other words, $F$
behaves from an homological viewpoint as a $p$-sheeted covering
from $\Sigma$ to $K$. Since $\CC^3$ is contractible, the class of
$\Sigma$ in $H_3(\CC^3)$ is zero. So the class $[K]$ in $H_3(\CC^2
-\{x\})$ is zero, hence a contradiction. Thus the main step in the
proof is the construction of the singular 3-chain $\Sigma$. We
proceed to this construction in sections {\bf{2-3-4}}.

Secondly, we are going to characterize from a topological
viewpoint the morphisms $F:\CC^n \rightarrow \CC^{n-1}$ which are
quotient maps of a $\CP$-action on $\CC^n$, and more generally
which are invariant with respect to such an action. Let $\varphi$
be an algebraic $\CP$-action in $\CC^n$ distinct from the
identity, i.e. $\varphi(t;x)\not\equiv x$. If $F:\CC^n \rightarrow
\CC^{n-1}$ is a dominant invariant morphism, then its generic
fibres are finite union of orbits of $\varphi$, hence their
connected components are contractible. Moreover if $F$ is a
quotient map for $\varphi$, then its generic fibres are connected
contractible, since they are one and only one orbit of $\varphi$.
Conversely polynomial maps with contractible generic fibres
correspond to $\CP$-actions. More precisely:

\begin{theorem} \label{cont}
Let $F:\CC^n \rightarrow \CC^{n-1}$ be a dominant polynomial map.
Assume that the connected components of its generic fibres are
contractible. Then there exists an algebraic $\CP$-action
$\varphi$, distinct from the identity, for which $F$ is invariant.
If moreover $F$ is non-singular in codimension 1 and its generic
fibres are connected, then $F$ is the quotient map of a
$\CP$-action on $\CC^n$.
\end{theorem}
Let $f_1,..,f_{n-1}$ be the coordinate functions of $F$. The main
idea is to introduce the following derivation: $$ \der : \CX
\rightarrow \CX, \quad R \mapsto J(R,f_1,..,f_{n-1})$$ where $J$
denotes the jacobian of $n$ functions in $n$ complex variables,
and to prove that $\der$ is locally nilpotent. Therefore its
integration leads to an algebraic $\CP$-action $\varphi$ on
$\CC^n$ for which $F$ is invariant, because $\der(f_i)=0$ for any
$i$. If the generic fibres are connected and $F$ is non-singular
in codimension 1, there remains to check that
$\CC[f_1,..,f_{n-1}]= \CX ^{\varphi}$, and this can be done by
using Zariski's Main Theorem. For more details, see section
{\bf{5}}.

As a consequence, theorem \ref{surj} can be rewritten in an
entirely topological way, as follows: If $F$ is a polynomial map
that is non-singular in codimension 1 and whose generic fibres are
connected contractible, then $F$ is surjective.

We end up this paper with two examples of polynomials maps which
are not surjective, and we will explain why in light of the
arguments given in the proof of theorem \ref{surj}. Moreover we
will see with the first example that there does not exist any
torical analogue of theorem \ref{cont}. More precisely there
exists a dominant map whose generic fibres are isomorphic to
$\CC^*$ and that is not invariant with respect to any
$\CC^*$-action. The second one is an example of a non-surjective
quotient map $F: \CC^4 \rightarrow \CC^3$. Both examples appear in
section {\bf{6}}.

\section{Some preliminary results}

We begin with some standard results concerning algebraic
$(\CC,+)$-actions on $\CC^3$. Recall that a $(\CC,+)$-action
$\varphi$ on $\CC^3$ induces a degree function $deg$ on
$\CC[x_1,x_2,x_3]$, defined by: $$ deg(R)=deg_t(R\circ
\varphi(t;x)) $$ Therefore $R$ is invariant with respect to
$\varphi$ if and only if $deg(R)=0$ or $R=0$. The existence of
this degree implies in particular that the ring of invariants is
factorially closed (see~\cite{Da}). This means that if a polynomial is
invariant, then all its irreducible factors are invariant. Let
$\Gamma$ be the complement in $\CC^2$ of $F(\CC^3)$. For any
polynomial $R$ in $\CC[x_1,..,x_n]$, we set by convention: $$
V(R)=\{x \in \CC^n, R(x)=0\},\quad D(R)=\{x \in \CC^n,
R(x)\not=0\} $$
\begin{lemma}
The set $\Gamma$ is at most finite.
\end{lemma}
\dem Since $f_1,f_2$ are algebraically independent, $F$ is a
dominant map. Moreover $\Gamma$ is a constructible set of
codimension $\geq 1$. Let us prove by absurd that $\Gamma$ has
codimension $\geq 2$. Suppose that $\Gamma$ contains a Zariski
open set $U$ of an irreducible curve in $\CC^2$. We may assume
that: $$U=D(Q)\cap V(P)$$ where $P$ is irreducible and $Q$ is not
divisible by $P$. Then $D(Q(F))\cap V(P(F))=\emptyset$, and
$V(P(F))\subset V(Q(F))$. By Hilbert's Nullstellensatz, there
exists an integer $n$ and a polynomial $R$ such that: $$Q(F)^n =
P(F)R$$ Since $\CC[f_1,f_2]$ is factorially closed, $R$ is of the
form $S(F)$, where $S$ is a polynomial. Therefore $$Q^n =P S$$ and
$P$ divides $Q$, hence a contradiction. \qed The following
lemma
is standard and asserts the existence of a rational slice for any
$(\CC,+)$-action on $\CC^3$ (see~\cite{Da},\cite{De},\cite{D-F}).

\begin{lemma} \label{fib}
Let $F$ be the quotient map of $\varphi$. Then there exists a hypersurface
$V(f)$ in $\CC^3$, and a principal open set $D(P)$ in $\CC^2$ such that
$F: V(f)\cap D(P(F))\rightarrow D(P)$ is an isomorphism.
\end{lemma}
\dem Let $\partial$ be the locally nilpotent derivation generating
$\varphi$. Since $\partial\not=0$, there exists a polynomial $f$
such that $\partial(f)\not=0$ and $\partial ^2(f)=0$. Since
$\CC[f_1,f_2]$ is the kernel of $\partial$, there exists a
polynomial $P$ such that $\partial(f)=P(F)$. By induction on the
degree, we easily check that every polynomial $R$ can be written
as $P(F)^n R = T(f,f_1,f_2)$, where $T$ is an element of
$\CC[x_1,x_2,x_3]$. This yields the equality: $$
\CC[x_1,x_2,x_3]_{P(F)} = \CC[f,f_1,f_2]_{P(F)} $$ So the map
$G=(f,f_1,f_2)$ defines an isomorphism from $D(P(F))$ to $\CC
\times D(P)$. Moreover $G$ maps $V(f)\cap D(P(F))$ on $\{0\}\times
D(P)$, and its restriction is equal to $F$ via the identification
$\{0\}\times D(P)\simeq D(P)$. Therefore $F: V(f)\cap
D(P(F))\rightarrow D(P)$ is an isomorphism. \qed

\section{Construction of coverings}

Let us denote by $Q$ an irreducible polynomial in $\CC[t_1,t_2]$,
and by $F$ a polynomial map from $\CC^3$ to $\CC^2$ that is
non-singular in codimension 1, that is whose singular set has
codimension $\geq 2$ in $\CC^3$. In this section, we will show how
to construct some coverings over a neighborhood of a compact set
contained in $V(Q)$. More precisely:

\begin{proposition} \label{sur1}
Let $F:\CC^3 \rightarrow \CC^2$ be a polynomial map that is nonsingular
in codimension 1. Let $Q$ be an irreducible polynomial in $\CC[t_1,t_2]$.
Then
there exists a Zariski open set $U$ of $V(Q)$ satisfying the following
property: For any compact
set $K$ contained in $U$, there exist an analytic subvariety $X_K$ in
$\CC^3$
and an open set $U_K$ in $\CC^2$, containing $K$, such that $F: X_K
\rightarrow U_K$
is a finite unramified covering.
\end{proposition}
This result applies in particular for any quotient map of an algebraic
$(\CC,+)$-action
on $\CC^3$, since Daigle proved in \cite{Da} that any such map is
non-singular in
codimension 1.

\begin{lemma} \label{sur2}
There exists a plane $H$ in $\CC^3$ and a point $x$ in
$H\cap V(Q(F))$ such that $F:H\rightarrow \CC^2$ is non-singular at $x$.
\end{lemma}
\dem Since $F$ is nonsingular in codimension 1, there exists a
point $x$
in the hypersurface $V(Q(F))$ such that $dF(x)$ has rank 2. In particular
the wedge product $df_1 \wedge df_2 (x)$ is non-zero. So there exists a
linear form $l$ on $\CC^3$ such that $dl\wedge df_1 \wedge df_2 (x)\not=0$.
Let us set:
$$
H=V(l-l(x))
$$
By construction $F:H\rightarrow \CC^2$ is non-singular at $x$, and $x$
belongs
to $H\cap V(Q(F))$.
\qed

\begin{lemma} \label{sur3}
Let $x$ and $H$ be a point and a plane in $\CC^3$ satisfying the
conditions of the previous lemma. Then there exists an irreducible curve
$C$ in $H$, passing through $x$ such
that the map $F:C \rightarrow V(Q)$ is dominant.
\end{lemma}
\dem Denote by $F_H$ the restriction map $F:H \rightarrow \CC^2$.
By assumption $F_H$
is smooth at the point $x$. So $F_H$ is dominant, $F^{-1}(V(Q))$ cannot
contain $H$ and
$F^{-1}(V(Q))\cap H$ is a union of irreducible curves. Since $Q(F(x))=0$,
there exists an irreducible curve passing through $x$ and contained in
$F^{-1}(V(Q))\cap H$. Let us fix such a curve and denote it by $C$.
Let us show by absurd that the restriction $F:C \rightarrow V(Q)$
is dominant. Assume it is not. Then $F$ maps $C$ to a point, and $F:C
\rightarrow V(Q)$
is everywhere singular. For any smooth point $x'$ of $C$, the differential
$dF(x')$ must vanish on $T_{x'} C$. So $dF(x')$ must have rank $<2$, and
the smooth part of $C$ is contained in the singular set $Sing(F_H)$. Since
this set is closed, $C$ is contained
in $Sing(F_H)$. But this is impossible because $x$ belongs
to $C$ and is not a singular point of $F_H$.
\qed

\begin{lemma} \label{sur4}
Let $x$ be a point, $H$ be a plane and $C$ be an irreducible curve
satisfying the
conditions of the previous lemmas. Then there exists a Zariski open set
$U$ of $V(Q)$ such that:
\begin{itemize}
\item{$F:F^{-1}(U)\cap C \rightarrow U$ is proper for the metric topology,}
\item{$U$ does not contain any critical value of $F:H\rightarrow \CC^2$.}
\end{itemize}
\end{lemma}
\dem The singular set of $F_H$ is closed in $H$ and does not
contain $x$. Since $x$ belongs to $C$ and $C$ is irreducible, the
intersection
$C\cap Sing(F_H)$ is at most finite. Let $U'$ be a Zariski
open set of $V(Q)$ such that $U'$ does not meet the finite set
$F(C\cap Sing(F_H))$. By assumption, $F:F^{-1}(U')\cap C \rightarrow U'$
is a dominant map
of irreducible curves. So this is a quasi-finite morphism, and there
exists a Zariski open set $U$ in $U'$ such that $F:F^{-1}(U)\cap C
\rightarrow U$ is finite, hence proper for the metric topology.
\qed \ \\
{\it {Proof of proposition \ref{sur1}:}} Let $F:\CC^3 \rightarrow \CC^2$
be a
polynomial map that is nonsingular in codimension 1. Let $Q$ be an
irreducible polynomial in $\CC[t_1,t_2]$. Let $H$ and $C$ be the plane in
$\CC^3$ and the irreducible curve found in the previous lemmas. Let $U$ be
the  Zariski open set of $V(Q)$ satisfying the conditions of lemma
\ref{sur4}.
Let $K$ be a compact set contained in $U$. Since $F:F^{-1}(U)\cap C
\rightarrow U$ is a proper map, $L=F^{-1}(K)\cap C$
is compact. Since $F^{-1}(U)$ does not meet the singular set of $F_H$,
there exists a relatively compact open set $U_1$ of $H$ that contains $L$
and does not meet $Sing(F_H)$. Then the restriction map:
$$
F:\overline{U_1} \rightarrow F(\overline{U_1})
$$
is proper because its source is compact. Moreover the set $F(U_1)$
contains $K$,
and is open because $F:U_1 \rightarrow \CC^2$ is non-singular. By the
localisation lemma (see~\cite{Ch}, p.29), there exist two open sets $X_K$
of $H$
and $U_K$ of $\CC^2$, containing $L$ and $K$ respectively, such that $X_K$
is contained in $U_1$ and the map $F:X_K \rightarrow U_K$ is proper for
the metric
topology. By construction $X_K$ is an analytic subvariety of $\CC^3$, and
the map $F:X_K
\rightarrow U_K$ is proper and non-singular. Since its fibres are compact
analytic
sets, they are finite (see~\cite{Ch}). Therefore $F:X_K \rightarrow U_K$
is a finite
unramified covering. \qed

\section{Proof of the first theorem}

From now on, we assume that the quotient map $F$ is not
surjective, or in other words that $\Gamma\not=\emptyset$. Up to a
translation, we may suppose that $\Gamma$ contains the origin in
$\CC^2$. In what follows, we will always consider singular
homology with integer coefficients.

Since $F$ is a continuous map from $\CC^3$ to $\CC^2 - \{0\}$, it
induces a morphism $F_*$ from the space of singular 3-chains in
$\CC^3$ to the space of singular 3-chains in $\CC^2 - \{0\}$. If
$\Delta^3$ denotes the standard 3-simplex, recall that a singular
3-chain ${\cal{K}}$ in a topological space $X$ is a formal sum:
$${\cal{K}}= \sum n_{\alpha} \Delta_{\alpha}$$ where each
$n_{\alpha}$ is an integer and each $\Delta_{\alpha}$ is a
continuous map from $\Delta^3$ to $X$. $\Delta_{\alpha}$ is called
a singular 3-simplex and its image is denoted by
$\Delta_{\alpha}(\Delta^3)$. In this section, we are going to
construct two singular 3-chains $\Sigma$ in $\CC^3$ and $K$ in
$\CC^2 - \{0\}$ such that:

\begin{itemize}
\item{The boundaries $\der \Sigma$ and $\der K$ are equal to zero, and
there exists an integer
$p>0$ such that $F_*(\Sigma)=pK$,}
\item{The class of $K$ in $H_3(\CC^2 - \{0\})\simeq \mathbb Z$ is a
generator of this group.}
\end{itemize}
Assume this is done for the moment. Since $\Sigma$ and $K$ have no
boundaries, they define homological classes $[\Sigma]$ and $[K]$
in $H_3(\CC^3)$ and in $H_3(\CC^2 - \{0\})$ respectively. Moreover
if $F_*$ denotes the morphism induced by $F$ on singular homology,
we get $F_*([\Sigma])=p[K]$. Since $\CC^3$ is contractible, we
have $[\Sigma]=0$ and $p[K]=0$, hence contradicting the fact that
$[K]$ is a generator of $H_3(\CC^2 - \{0\})\simeq \mathbb Z$.

So in order to complete the proof of theorem \ref{surj}, there
only remains to construct these singular chains. We proceed to
their construction in the following subsections.

\subsection{Construction of $K$}

Let $P$ be the polynomial appearing in lemma \ref{fib}, and let
$P_1,..,P_s$ be its irreducible factors in $\CC[t_1,t_2]$. For
each factor $P_i$, we denote by $U_i$ a Zariski open set of
$V(P_i)$ satisfying the conditions of proposition \ref{sur1}. Let
$S$ be a 3-sphere in $\CC^2$ centered at the origin in $\CC^2$.
Since the sets $V(P) - \cup_i U_i$, $\Gamma$ and $U_i\cap U_j$ for
$i\not=j$ are at most finite, we can choose its radius small
enough so that:
\begin{itemize}
\item{$S$ does not meet the sets $V(P) - \cup_i U_i$, $\Gamma$ and
$U_i\cap U_j$ for $i\not=j$,}
\end{itemize}
As every open set $U_i$ is connected for the metric topology, and
$U_i\cap U_j \cap S=\emptyset$ if $i\not=j$, every connected
component of $V(P)\cap S$ is contained in one and only one open
set $U_i$. Let us denote by $\gamma_1,..,\gamma_r$ the connected
components of $V(P)\cap S$. Since every $\gamma_i$ is contained in
an open set $U_j$, there exists an open set $U_{\gamma_i}$
containing $\gamma_i$ and satisfying the conditions of proposition
\ref{sur1}. Thus $S$ is covered by the open sets $D(P)$ and
$U_{\gamma_i}$, $1\leq i \leq r$.

\begin{lemma} \label{triang}
There exists a singular 3-cycle $K= \sum
n_{\alpha}\Delta_{\alpha}$ in $\CC^2 -\{0\}$ satisfying the
following conditions:
\begin{itemize}
\item{Every image $\Delta_{\alpha}(\Delta^3)$ is contained either in
$S\cap D(P)$ or
in one of the sets $S\cap U_{\gamma_i}$,}
\item{Every image $\Delta_{\alpha}(\Delta^3)$ cannot meet two different
sets $\gamma_i$
and $\gamma_j$,}
\item{If $\Delta_{\alpha}(\Delta^3)$ meets $\gamma_i$ and
$\Delta_{\beta}(\Delta^3)$ meets
$\gamma_j$ with $i\not=j$, then $\Delta_{\alpha}(\Delta^3)\cap
\Delta_{\beta}(\Delta^3)=\emptyset$,}
\item{The homological class of $K$ is a generator of $H_3(\CC^2 -\{0\})$.}
\end{itemize}
\end{lemma}
\dem Since the sphere $S$ is a deformation retract of $\CC^2
-\{0\}$, the inclusion map $i: S \hookrightarrow \CC^2 -\{0\}$
induces an isomorphism: $$ i_*: H_3(S) \rightarrow H_3(\CC^2
-\{0\})$$ So we can find a singular 3-cycle $K'$ generating
$H_3(\CC^2 -\{0\})$ of the form: $$K'= \sum
n'_{\alpha}\Delta'_{\alpha}$$where the image of every
$\Delta'_{\alpha}$ lies in $S$. Let $d$ be the distance function
defined by the canonical Hermitian metric on $\CC^2$. We provide
$S$ with the metric topology induced by the embedding
$i:S\hookrightarrow \CC^2$. Since $S$ is compact and covered by
the open sets $D(P),U_{\gamma_1},..,U_{\gamma_r}$, there exists an
$\epsilon$ such that any ball of radius $\leq \epsilon$ in $S$ is
contained in one of these open sets. For any compact sets
$\gamma,\gamma'$ in $\CC^2$, we denote by $dist(\gamma, \gamma')$
the distance between these two sets. Up to choosing a smaller
$\epsilon$, we may even assume that: $$\epsilon \leq
\frac{dist(\gamma_i,\gamma_j)}{3}$$ whenever $i\not=j$. By
performing enough barycentric subdivisions of every simplex
$\Delta'_{\alpha}$ in $K'$, we can get a new 3-cycle $K$,
homologous to $K'$, such that:

$$K= \sum n_{\alpha}\Delta_{\alpha}$$ where every image
$\Delta_{\alpha}(\Delta^3)$ is contained in $S$ and has diameter
$\leq \epsilon$. By construction, the homological class of $K$ is
a generator of $H_3(\CC^2 -\{0\})$. Moreover since every set
$\Delta_{\alpha}(\Delta^3)$ has diameter $\leq \epsilon$, it is
contained in a ball of radius $\epsilon$. Hence
$\Delta_{\alpha}(\Delta^3)$ is contained in one of the open sets
$S\cap D(P),S\cap U_{\gamma_1}..,S\cap U_{\gamma_r}$ in $S$. The
other two conditions are as easy to check. \qed Let $K$ be the
singular 3-chain of the previous lemma. By the second condition of
this lemma, we can perform the following partition:
\begin{itemize}
\item{$\{\Delta_{(i,j)}\}_j$ is the set of 3-simplices of $K$
meeting $\gamma_i$,}
\item{$\{\Delta_k\}_k$ is the set of 3-simplices of $K$
meeting none of the $\gamma_i$.}
\end{itemize}
That enables us to rewrite this singular 3-cycle in the following
way: $$ K= \sum_{i,j} n_{(i,j)} \Delta_{(i,j)} + \sum_{k}
n_k\Delta_{k} $$ Note that by the third condition of the lemma,
the images of $\Delta_{(i,j)}$ and $\Delta_{(i',j')}$ intersect
only if $i=i'$.

\subsection{Construction of $\Sigma$}

In this subsection, we construct the singular 3-chain $\Sigma$ by
lifting the 3-simplices of $K$ in a suitable way. Let
$X_{\gamma_i}$ be the analytic variety given by proposition
\ref{sur1}, and let $p_i$ be the degree of the unramified
covering: $$ F: X_{\gamma_i}\rightarrow U_{\gamma_i} $$ Since the
image of every $\Delta_{(i,j)}$ is contained in $U_{\gamma_i}$, we
can lift it in $p_i$ different ways. More precisely there exist
$p_i$ different maps $\Delta^l _{(i,j)}:\Delta^3\rightarrow
X_{\gamma_i}$ making the following diagram commute:

$$\ \ \ \ \ \ \ \ \ \ \ \ X_{\gamma_i}$$ $$ \ \ \ \ \ \ \nearrow \
\ \ \ \ \big\downarrow F$$ $$\Delta^3\ \longrightarrow \
U_{\gamma_i}$$ where the arrow at the bottom stands for the map
$\Delta_{(i,j)}$. Let us denote by $\Delta ^1 _k$ the lifting of
the 3-simplex $\Delta_k$ via the isomorphism: $$ F:V(f)\cap
D(P(F))\rightarrow D(P) $$ More precisely, if $G$ is the
restriction of $F$ to $V(f)\cap D(P(F))$, then $\Delta^1 _k$ is
the map $G^{-1}\circ \Delta_k$. This yields the other
following commutative diagram:

$$\ \ \ \ \ \ \ \ \ \ \ \ V(f)\cap D(P(F))$$ $$ \ \ \ \ \ \
\nearrow \ \ \ \ \ \big\downarrow F$$ $$\Delta^3 \ \longrightarrow
\ D(P)$$ where the arrow at the bottom stands for the map
$\Delta_k$. Now, and {\em this is the key-point of the
construction}, we are going to modify these simplices so that the
boundary of $\cup_{i,j,l}\Delta^l _{(i,j)}$ coincides with the
boundary of $\cup_{k}\Delta^1 _{k}$. That will enable us to get a
singular 3-chain $\Sigma$ {\em with no boundary}. In order to do
so, we will use the fact that the generic fiber of $F$ is
contractible, in the following way. Let $V$ be the complement in
$S$ of $\cup_k \Delta_{k}(\Delta^3)$. By construction $V$ is an
open neighborhood of the union $\cup_i \gamma_{i}$ in $S$. There
exists a continuous function $g$ on $S$ such that:
\begin{itemize}
\item{$g$ is equal to 1 outside $V$,}
\item{$g$ vanishes in a neighborhood of each $\gamma_i$.}
\end{itemize}
If $\varphi$ is the algebraic $(\CC,+)$-action of the beginning,
we define the map $L$ on $F^{-1}(S)$ by the formula: $$
L(x)=\varphi(t(x);x),\quad t(x)=\frac{-f(x)g(F(x))}{P(F(x))} $$
outside $\cup_i F^{-1}(\gamma_i)$, and $L$ is the identity on
$\cup_i F^{-1}(\gamma_i)$. Note that $L$ is continuous since $g$
vanishes on a neighborhood of each $\gamma_i$ in $S$. Moreover
$F\circ L=F$ on $S$. The new 3-simplices ${\cal{D}}^l _{(i,j)}$
and ${\cal{D}}^1 _{k}$ are given by the formulas: $$ {\cal{D}}^l
_{i,j} =L\circ \Delta^l _{(i,j)} \quad \mbox{and} \quad
{\cal{D}}^1 _{k}=L\circ \Delta^1 _{k} $$ By construction we get:
$$F\circ {\cal{D}}^l _{(i,j)}= F\circ \Delta^l
_{(i,j)}=\Delta_{(i,j)} \quad \mbox{and} \quad F\circ {\cal{D}}^1
_{k}= F\circ \Delta^l _{k}= \Delta _k $$ We set $p=\prod_i p_i$
and define the singular 3-chain $\Sigma$ by the sum: $$ \Sigma =
\sum_{i,j,l} \frac{p}{p_i} n_{(i,j)}{\cal{D}}^l _{(i,j)} + p
\sum_k n_k{\cal{D}}^1 _{k} $$

\subsection{Properties of these singular 3-chains}

We are going to derive the properties announced at the beginning
of this section, and then conclude the proof of theorem
\ref{surj}. Recall that a face of a 3-simplex $\Delta$ is the
restriction of $\Delta$ to one of the faces of $\Delta^3$. By
extension a face of a 3-chain is a face of one of the 3-simplices
of its decomposition. Let $\Sigma'$ be the singular 3-chain: $$
\Sigma' = \sum_{i,j,l} \frac{p}{p_i} n_{(i,j)}\Delta ^l _{(i,j)} +
p \sum_k n_k \Delta^1 _{k} $$ For commodity we introduce the
following sets:
\begin{itemize}
\item{$E$ is the set of faces $\delta$ of $\Sigma$ such that $F(\delta)$
belongs to
the boundary of both a $\Delta_{(i,j)}$ and a $\Delta_k$,}
\item{$E'$ is the set of faces $\delta'$ of $\Sigma'$ such that
$F(\delta')$ belongs to
the boundary of both a $\Delta_{(i,j)}$ and a $\Delta_k$.}
\end{itemize}

\begin{proposition} \label{sb4}
If $F_{*}$ is the morphism induced by $F$ on the space of singular
3-chains, then $F_{*}(\Sigma)= p K$.
\end{proposition}
\dem By construction we have the following relations: $$
F_{*}({\cal{D}}^l _{(i,j)})= \Delta_{(i,j)} \quad \mbox{and} \quad
F_{*}({\cal{D}}^1 _{k})= \Delta_{k} $$ Since every 3-simplex
$\Delta_{(i,j)}$ has been lifted $p_i$ times, and every 3-simplex
$\Delta_{k}$ has been lifted once, we obtain:
\begin{equation*}
\begin{split}
F_{*}(\Sigma)& = \sum_{i,j} \frac{p}{p_i}n_{(i,j)}\left(\sum_l
F_{*}({\cal{D}}^l _{(i,j)})\right ) + p\sum_{k}n_k
F_{*}({\cal{D}}^1 _{k}) \\ & = \sum_{i,j} \frac{p}{p_i} p_i
n_{(i,j)} \Delta _{(i,j)} + p\sum_{k} n_k \Delta_{k}  \\ & =  p
\sum_{i,j} n_{(i,j)} \Delta _{(i,j)} + p\sum_{k} n_k \Delta_{k} \\
& = pK
\\
\end{split}
\end{equation*}
\qed

\begin{lemma} \label{sb3}
Let $\delta_1,\delta_2$ be any 2-faces of $\Sigma$ such that
$\delta=F(\delta_1)=F(\delta_2)$ belongs to the boundary of a
$\Delta_{k}$. Then $\delta_1=\delta_2$.
\end{lemma}
\dem Let us show that $\delta_1$ is equal to the map $G^{-1}\circ
\delta$ (see the previous subsection). If $\delta_1$ belongs to
the boundary of a ${\cal{D}}^1_{k}$, then $\delta_1=G^{-1}\circ
\delta$ by construction. If $\delta_1$ belongs to the boundary of
a ${\cal{D}}^l _{(i,j)}$, then $\delta_1$ is of the form
$L\circ\delta' _1$, where $\delta' _1$ is a face of $\Delta ^l
_{i,j}$. Since $F(\delta'_1)=F(\delta_1)=\delta$ belongs to the
boundary of a $\Delta_{k}$, the function $g$ is equal to 1 on the
image of $F(\delta'_1)$. So $g\circ F \circ \delta'_1 =1$ and we
have the following equality for any point $x$ in the image of
$\delta'_1$: $$ L(x)=\varphi(-f(x)/P(F(x));x) $$ By using the
exponential map, we get: $$f\circ \varphi(t;x)= f(x) + P(F(x))t$$
After substitution, that implies: $$f(L(x))= 0$$ Therefore $f\circ
\delta_1=f\circ L \circ \delta'_1=0$ and the image of $\delta_1$
lies in the set $V(f)$. Since $\delta$ is a face of a $\Delta_k$,
and $\Delta_k$ does not meet the hypersurface $V(P)$, the function
$P\circ \delta$ never vanishes. Since $F(\delta_1)=\delta$, the
function $P\circ F \circ \delta_1$ never vanishes and the image of
$\delta_1$ lies in the intersection $V(f)\cap D(P(F))$. Since
$F\circ \delta_1=\delta$, and since the map $G$ defined as the
restriction: $$F:V(f)\cap D(P(F)) \rightarrow D(P)$$ is an
isomorphism, $\delta_1$ is equal to $G^{-1}\circ \delta$.
\qed

\begin{lemma} \label{sb1}
For any $i$, the boundary of $\sum_{j,l} n_{(i,j)}{\cal{D}}^l
_{(i,j)}$ is a linear combination of elements of $E$.
\end{lemma}
\dem Assume first that the boundary of $\sum_{j,l} n_{(i,j)}\Delta
^l _{i,j}$ can be written as: $$
\partial \left (\sum_{j,l} n_{(i,j)}\Delta ^l _{i,j}\right )=
\sum_{\delta'\in E'}
n_{\delta'} \delta'$$ Let $L_{*}$ be the morphism induced by $L$.
Since $L \circ \Delta ^l _{i,j}={\cal{D}}^l _{(i,j)}$ and $L \circ
\Delta ^1 _{k}={\cal{D}}^1 _{k}$, $L_{*}$ maps elements of $E'$ to
elements of $E$. Since the boundary operator commutes with
$L_{*}$, that implies: $$
\partial \left (\sum_{j,l} n_{(i,j)}{\cal{D}} ^l _{i,j}\right )= \sum_{\delta\in E}
\left (\sum_{L_{*}(\delta')=\delta}n_{\delta'}\right ) \delta $$
So there only remains to show that the boundary of $\sum_{j,l}
n_{(i,j)}\Delta ^l _{i,j}$ is a linear combination of elements of
$E'$. Let us prove that any face $\delta'$ not belonging to $E'$
cannot appear with a non-zero coefficient into that boundary.

Let $\delta'$ be a face of a $\Delta ^l _{i,j}$, that does not
belong to $E'$. Then $\delta=F(\delta')$ is a face of
$\Delta_{(i,j)}$. Moreover $\delta$ is not a face of any
$\Delta_{k}$. By lemma \ref{triang}, the only simplices of
$\Sigma'$ that may have $\delta$ as a face are of the form
$\Delta_{i,j'}$. We write them as
$\Delta_{i,j_1},..,\Delta_{i,j_r}$. We now use the lifting
property of the covering: $$F:X_{\gamma_i} \rightarrow U_i$$ For
any $\alpha$, there exists a unique lifting $\Delta^{l_{\alpha}}
_{i,j_{\alpha}}$ of $\Delta_{i,j_{\alpha}}$ such that $\delta'$ is
one of its faces. So $\Delta^{l_1} _{i,j_1},..,\Delta^{l_r}
_{i,j_r}$ are the only simplices of $\Sigma'$ and of $\sum_{j,l}
n_{(i,j)} \Delta ^l _{i,j}$ having $\delta'$ as a face. Let
$\epsilon_{\alpha}$ be the coefficient of $\delta$ in the boundary
of $\Delta _{i,j_{\alpha}}$. Since the singular 3-chain $K$ has no
boundary, we have: $$
 \sum_{\alpha} n_{(i,j_{\alpha})}\epsilon_{\alpha}=0
$$ But $\epsilon_{\alpha}$ is also the coefficient of $\delta'$ in
the boundary of $\Delta^{l_{\alpha}} _{i,j_{\alpha}}$. Therefore,
the coefficient of $\delta'$ in the boundary of $\sum_{j,l}
n_{(i,j)}\Delta ^l _{i,j}$ is equal to the number given above,
hence zero, and the result follows. \qed

\begin{lemma} \label{sb2}
The boundary of $\sum_k n_k {\cal{D}}^1 _{k}$ is a linear
combination of elements of $E$.
\end{lemma}
\dem The proof is entirely similar to the proof of the previous
lemma. The only difference is the use of the isomorphism
$F:V(f)\cap D(P(F))\rightarrow D(P)$ in place of the covering
$F:X_{\gamma_i}\rightarrow U_i$ for the definition of the
$\Delta^1 _k$. \qed

\begin{proposition}
The singular 3-chain $\Sigma$ has no boundary.
\end{proposition}
\dem By applying lemmas \ref{sb1} and \ref{sb2} to the definition
of $\Sigma$, we see that the boundary of $\Sigma$ can be written as:
$$
\partial \Sigma = \sum_{\delta' \in E} n_{\delta'} \delta'
$$ Since the boundary operator commutes with the morphism $F_{*}$
induced by $F$, we get by lemma \ref{sb4}: $$ \sum_{\delta \in
F(E)} \left (\sum_{F_{*}(\delta')=\delta}n_{\delta'}\right )
\delta= p\partial(K)=0 $$ Thus all the sums
$\sum_{F_{*}(\delta')=\delta}n_{\delta'}$ are equal to zero. By
lemma \ref{sb3}, for any face $\delta$ in $F(E)$, there exists a
unique face $\delta'$ in $E$ such that $F(\delta')=\delta$. This
implies the equality $n_{\delta'}= 0$ for any $\delta'$, and the
result follows. \qed

\section{Morphisms with contractible generic fibres}

In this section, we pass on to polynomial maps with contractible
fibres, and we are going to prove theorem \ref{cont}. We begin
with the following lemma, which corresponds to the first assertion
of this theorem.

\begin{lemma}
Let $F$ be a dominant polynomial map from $\CC^{n}$ to
$\CC^{n-1}$. Assume that the connected components of its generic
fibres are contractible. Then there exists an algebraic
$\CP$-action $\varphi$, distinct from the identity, for which $F$
is invariant. In particular, the generic fibres of $F$ are finite
unions of orbits of $\varphi$.
\end{lemma}
\dem Write $F=(f_1,..,f_{n-1})$, and let $\der$ be the derivation
on $\CX$ defined for any $R$ by: $$ \der(R) = J(R,f_1,..,f_{n-1})
$$ where $J$ denotes the jacobian of $n$ functions in $n$
variables. If we show that $\der$ is a locally nilpotent
derivation, then it will generate an algebraic $\CP$-action
$\varphi$ on $\CC^n$ for which each $f_i$ is invariant, because
$\der(f_i)=0$ for any $i$. Moreover $\varphi$ will be distinct
from the identity because $\der\not=0$. So let us prove by absurd
that $\der$ is locally nilpotent.

Assume there exists a polynomial $R$ such that $\der ^k (R)\not=0$
for any $k$. By assumption on $F$, there exists a Zariski open set
$U$ in $\CC^{n-1}$ such that:
\begin{itemize}
\item{For any $y$ in $U$, $y$ is not a critical value of $F$ and $F^{-1}(y)$ is not empty,}
\item{For any $y$ in $U$, the connected components of $F^{-1}(y)$
are contractible.}
\end{itemize}
For any $k$, let us set $U_k = F^{-1}(U) \cap D(\der ^k(R))$. Then
$U_k$ is a non-empty Zariski open set which is dense in $\CC^n$
for the metric topology. By Baire's property of complete
topological spaces, we get: $$ \cap _{k \geq 0} \; U_k \not=
\emptyset $$ Let $x$ be a point of this intersection, and let $C$
be the connected component of $F^{-1}(F(x))$ containing $x$. Since
$\der(f_i)=0$ for any $i$, $\der$ corresponds to a vector field
that is tangent to $F^{-1}(F(x))$, hence to $C$. Therefore it
induces a derivation $\Delta$ on the ring $\CC[C]$ of regular
functions on $C$. Since $C$ is a smooth contractible algebraic
curve, it is isomorphic to $\CC$, and $\Delta$ appears as a
derivation on $\CC[t]$. Write it as: $$\Delta=
P(t)\frac{\der}{\der t} $$ By construction the singularities of
$\der$, considered as a vector field, are the singular points of
$F$. Since $F(x)$ is not a critical value of $F$, $\der$ has no
singularities along $C$. Thus the polynomial $P(t)$ must never
vanish, hence it is a constant. Therefore $\Delta$ is locally
nilpotent on $\CC[t]$, and there exists an order $k$ such that: $$
\Delta ^k (R) = \der ^k (R)_{|C} = 0 $$ In particular $\der^k
(R)(x)=0$, hence contradicting the construction of $x$. \qed

\noindent {\it {Proof of theorem \ref{cont}}}: There only remains
to show the second assertion of this theorem. Let $F$ be a
dominant polynomial map that is non-singular in codimension 1, and
whose generic fibres are connected contractible. By the previous
lemma, there exists a $\CP$-action $\varphi$ distinct from the
identity and for which $F$ is invariant. Let us check that: $$\CX
^{\varphi}= \CC[F]$$ Let $R$ be an invariant polynomial. Since the
generic fibres of $F$ are smooth and connected, each of them is
exactly one and only one orbit of $\varphi$. Let $U$ be a Zariski
open set in $\CC^{n-1}$ such that, for any $y$ in $U$, $F^{-1}(y)$
is an orbit of $\varphi$. Since $R$ is invariant, it is constant
along any such fibre. Consider the following correspondence:
$$\alpha: U \longrightarrow \CC, \quad y \longmapsto \mbox{"value
of $R$ along $F^{-1}(y)$"} $$ Its graph corresponds to the image
of $F^{-1}(U)$ by the map $(f_1,..,f_{n-1},R)$. This is a
constructible set whose Zariski closure is irreducible. Thus
$\alpha$ defines a rational correspondence in the sense of
Zariski. By Zariski's Main Theorem (see \cite{Mum}), $\alpha$ is
rational and $R$ can be written as $R=\alpha(F)$. Let us prove by
absurd that $\alpha$ is a polynomial.

Assume that $\alpha= A/B$, where $A$ and $B$ have no common
factors and $B$ is not constant. Since $R=A(F)/B(F)$ is a
polynomial, $A(F)$ and $B(F)$ have a common irreducible factor
$h$. So $F$ maps the hypersurface $V(h)$ into the set $V(A)\cap
V(B)$. For any smooth point $x$ of $V(h)$, we get: $$ rank
(d\{F_{|V(h)} \}(x))= rank(dF(x) _{|T_x V(h)})\leq dim\; V(A)\cap
V(B) \leq (n-3) $$ Since $T_x V(h)$ is an hyperplane of $\CC^n$,
this yields: $$ rank(dF(x))\leq (n-2) $$ By upper semi-continuity,
$dF(x)$ has rank $\leq (n-2)$ for any point $x$ of $V(h)$.
Therefore the singular set of $F$ contains the hypersurface
$V(h)$, hence contradicting the fact that $F$ is nonsingular in
codimension 1. \qed

\section{Two examples}

Finally we are giving two examples of polynomial maps between
affine spaces which are not surjective, and we are trying to
understand why. Moreover we are going to see with the first
example that there is no torical analogue of theorem \ref{cont},
namely that a polynomial map $F:\CC^n \rightarrow \CC^{n-1}$,
whose generic fibres are isomorphic to $\CC ^*$, need not be the
quotient map of an algebraic $\CC^*$-action on $\CC^n$.

\subsection{First example}

The construction of the singular 3-chains $K$ and $\Sigma$, which
is the main argument of the proof of theorem \ref{surj}, is made
possible because first $F$ is non-singular in codimension 1 and
second, the generic fibre of $F$ is contractible. This is clear by
theorem \ref{cont}. These are the reasons why we can lift the
singular 3-chain $K$, and then ajust the boundaries of the
singular 3-simplices forming $\Sigma$ so as to assure that
$\Sigma$ has no boundary. Consider the following map: $$F: \CC^3
\rightarrow \CC^2, \quad (x,y,z) \rightarrow (1+xz,y+z+xyz)$$ Then
$F$ is not surjective because its image is the set $\CC^2 -\{0\}$.
Moreover its singular set is the line $\{x=0, z=0\}$ in $\CC^3$,
and so $F$ is non-singular in codimension 1. The obstruction to
surjectivity lies in the fact that the generic fibre of $F$ is
isomorphic to $\CC^*$, hence it is not contractible.

As in theorem \ref{cont}, we might expect that $F$ is invariant
with respect to an algebraic $\CC^*$-action, that is a regular map
$\psi : \CC^* \times \CC^3 \rightarrow \CC^3$ such that $\psi(t;
\psi (s; p))= \psi (ts; p)$ for any $t,s,p$. We are going to see
that this is not the case. Assume that $F$ is invariant with
respect to such an action $\psi$, whose parameter in $\CC^*$ is
denoted by $t$. Then the polynomial $xz$ is invariant with respect
to $\psi$, and there exists an integer $r$ such that: $$ x \circ
\psi = t^r x, \quad z \circ \psi = t^{-r} z $$ Since $y+z+xyz$ is
invariant, this yields the equality: $$ (y\circ \psi - y)(1+xz) =
z (1 - t^{-r}) $$ and this is impossible since $(1+xz)$ cannot
divide $z(1-t^{-r})$ in $\CC[t,1/t,x,y,z]$.

\subsection{Second example}

Considering the class of algebraic $(\CC,+)$-actions on $\CC^n$
satisfying condition $(H)$ (see the introduction), we may ask if
theorem \ref{surj} extends with no restriction to higher
dimension, that is for $n>3$. More precisely, if $\varphi$ is an
algebraic $(\CC,+)$-action on $\CC^n$ satisfying condition $(H)$,
is its quotient map always surjective ?

The answer is no. Let us denote by $x,y,u,v$ the coordinates in
$\CC^4$, let $t$ be a parameter in $\CC$ and consider the
$(\CC,+)$-action on $\CC^4$ defined as follows: $$ \varphi
(t;x,y,u,v)=(x,y,u -t y, v+tx)$$ It is easy to check that its ring
of invariants is generated by the polynomials $x,y,xu+yv$. So
$\varphi$ satisfies condition $(H)$, and its quotient map is given
by: $$F: \CC^4 \longrightarrow \CC^3, \quad (x,y,u,v)\longmapsto
(x,y,xu+yv)$$ The map $F$ is not surjective, since its image is
the set: $$F(\CC^4)= \CC^3 - \{(x_1,x_2,x_3), x_1=x_2=0,
x_3\not=0\}$$ We may wonder why this map is not surjective, since
surjectivity is automatically satisfied for quotient maps if
$n=3$. In fact, given a singular 3-cycle $K$ in $F(\CC^4)$, we can
reproduce in exactly the same way our previous construction, and
find a singular 3-cycle $\Sigma$ in $\CC^4$ and an integer $p>0$
such that: $$F_*(\Sigma)=pK$$ But this will not lead us to any
contradiction as in the proof of theorem \ref{surj}, because the
group $H_3(F(\CC^4))$ is reduced to zero. Indeed the set
$F(\CC^4)$ is contractible, since it retracts by deformation to
the origin via the following map: $$ R: [0,1]\times F(\CC^4)
\longrightarrow F(\CC^4), \quad (t,x_1,x_2,x_3)\longmapsto
(tx_1,tx_2,tx_3)$$

\end{document}